\documentclass[12pt]{article}
\usepackage[psamsfonts]{amssymb}
\usepackage{amsthm}

\theoremstyle{plain}

\newtheorem{proposition}{Proposition}[section]
\newtheorem{corollary}{Corollary}[section]

\theoremstyle{definition}
\newtheorem{definition}{Definition}[section]
\newtheorem{example}{Example}[section]

\theoremstyle{remark}
\newtheorem*{remark}{\indent Remark}

\newcommand{\R}{\mathbf{R}}                     
\newcommand{\dd}{{\rm d}}

\title{Derivative relationships between volume and surface area of compact regions in $\mathbf{R}^d$}

\author{Jean-Luc Marichal\footnote{Corresponding author.}
\and Michael Dorff}

\date{}

\begin{document}
\maketitle

\begin{abstract}
We explore the idea that the derivative of the volume, $V$, of a region in $\R^d$ with respect to $r$ equals its surface area, $A$, where $r=d\,
\frac VA$. We show that the families of regions for which this formula for $r$ is valid, which we call homogeneous families, include all the
families of similar regions. We determine equivalent conditions for a family to be homogeneous, provide examples of homogeneous families made up
of non-similar regions, and offer a geometric interpretation of $r$ in a few cases.
\end{abstract}

\bigskip


\section{Introduction}

It is well known that there exists a remarkable derivative relationship between the area $A$ and the perimeter $P$ of a circle, namely $$\frac
{\dd A}{\dd r}=P,$$ where the variable $r$ represents the radius of the circle. It is natural to wonder whether such a derivative relationship
remains valid for other familiar shapes. At first glance, though, it does not even hold for the square when $r$ represents the side length.
However, it holds when $r$ represents half of the side length, that is, the radius of the inscribed circle.

In a similar manner, the derivative of the volume function of a sphere is equal to the surface area, that is,
$$\frac {\dd V}{\dd r}=A$$ and this relationship still holds for cubes if $r$ represents the radius of the inscribed sphere.

We show that by choosing an appropriate variable to calculate volume and area, namely
\begin{equation}\label{eq:tongintro}
r=d\,\frac{V}{A}
\end{equation}
(as recently suggested by Tong~\cite{To97}), we can generalize the derivative relationship to many compact regions in $\R^d$ $(d\geqslant 2)$.

Notice that, when we consider the derivative relationship of a given compact region, we actually consider a one-parameter family of similar
compact regions. For example, the derivative relationship for a sphere involves considering a sphere that grows in radius, that is, a family of
spheres.

Also, we can investigate families of non-similar regions. For example, consider a right circular cone in $\R^3$ whose base radius and height are
functions of a certain parameter $s$. We can calculate the volume $V(s)$ and the surface area $A(s)$ as functions of $s$ and then search for an
appropriate change of variable $r(s)$ for which the derivative relationship holds.

In this general case of possibly non-similar regions, we show that the derivative relationship always holds for the change of variable
\begin{equation}\label{eq:chvarr}
r(s)=\int\frac{V'(s)}{A(s)}\, \dd s.
\end{equation}

In this paper we mainly investigate one-parameter families of regions for which the change of variable reduces to (\ref{eq:tongintro}). We call
these families {\it homogeneous} families and we show that a family is homogeneous if and only if its regions have the same isoperimetric ratio.
In particular, any family of similar regions is homogeneous. We also show how to construct homogeneous families made up of non-similar regions.

The outline of this paper is as follows. In the next section we derive the change of variable formula (\ref{eq:chvarr}). In Section 3 we provide
characterizations of the class of homogeneous families. In Section 4 we show how to construct such families. In Section 5 we yield a geometric
interpretation of the variable (\ref{eq:tongintro}) for certain homogeneous families such as families of star-like polyhedra. Finally, in
Section 6 we provide Bonnesen-type isoperimetric inequalities constructed from this latter variable.

Surprisingly, derivative relationships between volume and area of compact regions have not been widely investigated. To our knowledge, only a
few researchers have worked on this interesting topic; see \cite{DoHa}, \cite{EmNe97}, \cite{Mi78}, \cite{St90}, \cite{To97}.

Throughout, we will use the notation $\R_+$ for the interval $(0,+\infty)$.

\section{Derivative relationship: the general case}

Let $d\geqslant 2$ be an integer. Consider a one-parameter family of compact regions in $\R^d$ with boundaries of finite measures,
$$\mathcal{F}:=\{R(s)\subset\R^d\mid s\in E\},$$ where $E$ is an open interval of the real line. We assume that with this
family is associated a strictly monotone and differentiable function $V:E\to\R_+$ and a continuous function $A:E\to\R_+$ such that, for any
$s\in E$, the values $V(s)$ and $A(s)$ represent respectively the volume and the surface area of region $R(s)$.

For the sake of convenience, such a family will be called a {\em smooth}\/ family.

Note that for plane figures in $\R^2$, we replace the volume $V(s)$ and the area $A(s)$ with the area $A(s)$ and the perimeter $P(s)$,
respectively.

The parameter $s$ can represent either a linear dimension, or an angle, or may have no geometric meaning.

\begin{example}\label{ex:cube}
Consider a (smooth) family of cubes in $\R^3$, with edge length $s\in\R_+$. In that case the volume and area functions are clearly given by
$V(s)=s^3$ and $A(s)=6s^2$, respectively. Of course, we could as well choose any positive function $\phi(s)$ of the parameter $s$ to represent
the edge length, thus leading to the new functions $V(s)=\phi(s)^3$ and $A(s)=6\phi(s)^2$. In such an alternative representation the parameter
$s$ may have no geometric interpretation.
\end{example}

Let us now show that, for any smooth family, it is always possible to find an appropriate variable of differentiation leading to the derivative
relationship between volume and surface area.

\begin{proposition}\label{prop:r}
Let ${\cal F}$ be a one-parameter smooth family of compact regions. Then there is a differentiable change of variable $r(s):E\to r(E)$, defined
as
\begin{equation}\label{eq:r}
r(s):=\int\frac{V'(s)}{A(s)}\, \dd s \qquad (s\in E)
\end{equation}
and unique within an additive constant, such that
\begin{equation}\label{eq:derrel}
\frac \dd{\dd r}\, V[s(r)]=A[s(r)] \qquad (r\in r(E)).
\end{equation}
\end{proposition}

\begin{proof}
The sign of the derivative
\begin{equation}\label{eq:rp}
r'(s) = \frac{V'(s)}{A(s)} \qquad (s\in E)
\end{equation}
is constant and $r(s)$ is a differentiable change of variable from $E$ to $r(E)$.

By the chain rule, we then have $$\frac \dd{\dd r}\, V[s(r)] = V'[s(r)]\, s'(r)=\frac{V'[s(r)]}{r'[s(r)]}=A[s(r)]$$ for all $r\in r(E)$. The
uniqueness of $r(s)$ follows immediately from the latter equality.
\end{proof}

From Eq.~(\ref{eq:r}) we immediately see that the variable of differentiation $r$ represents a linear dimension. Moreover, if $V(s)$ and $A(s)$
are replaced with
$$V_{\phi}(s)=V[\phi(s)]\qquad\mbox{and}\qquad
A_{\phi}(s)=A[\phi(s)],$$ respectively, where $\phi:E\to \phi(E)\subseteq E$ is a differentiable change of variable, then $r(s)$ is simply
replaced with
\begin{eqnarray*}
r_{\phi}(s) &=& \int\frac{V_{\phi}'(s)}{A_{\phi}(s)}\, \dd s = \int\frac{V'[\phi(s)]\, \phi'(s)}{A[\phi(s)]}\, \dd s = \int\frac{V'(t)}{A(t)}\,
\dd t\Big|_{t=\phi(s)}\\ &=& r[\phi(s)],
\end{eqnarray*}
which clearly shows that the change of variable $r(s)$ remains stable under any change of representation.

In Example~\ref{ex:cube}, with the family of cubes of edge lengths $s$, we have
$$r(s) = \frac s2+C,$$ for a constant $C\in\R$. When $C=0$, the variable $r$
represents the radius of the inscribed sphere. When $C\neq 0$, this radius is given by $r-C$. We then have
$$V[s(r)]=8r^3 \quad \mbox{and}\quad A[s(r)]=24r^2,$$ thus retrieving
Eq.~(\ref{eq:derrel}) with $E=r(E)=\R_+$.

Although the new variable $r$ represents a length, a geometric interpretation of it is not always immediate, as the following example shows.

\begin{example}\label{ex:rectfixw}
Consider a family of rectangles with fixed length $a>0$ and variable width $s>0$. Then we have $A(s)=as$, $P(s)=2s+2a$, and
$$r(s)=\frac a2\,\ln(2s+2a)+C.$$ In this case, no interpretation is known for the variable $r$.
\end{example}

As we will see in the subsequent sections, when the regions of ${\cal F}$ are all similar, $r$ takes a simpler form and can sometimes be
interpreted.

\begin{example}\label{ex:rectsim}
Consider a family of similar rectangles with length $s>0$ and width $ks$, where $k\in (0,1)$ is a fixed constant. Then we have $A(s)=ks^2$,
$P(s)=2s+2ks$, and
$$r(s)=\frac k{k+1}\, s+C.$$ In this case, $r$ is
the harmonic mean of the half-length and the half-width, i.e.,
$$
r(s)=H\big(\frac s2, k\,\frac s2\big)+C.
$$
\end{example}

Notice also that it is necessary that $V(s)$ be strictly monotone in $E$ for $r(s)$ to be a change of variable. In situations where $V(s)$ is
not strictly monotone in its domain, it is necessary to partition this domain into open subintervals $E$ in which $V(s)$ is strictly monotone.

\begin{example}\label{ex:vstrmon}
Consider a family of rhombi in $\R^2$ with sides of fixed length $a>0$ and a diagonal of variable length $s\in (0,2a)$. The perimeter $P(s)=4a$
is constant while the area
$$A(s)=s\, \sqrt{a^2-\frac{s^2}4}$$ is strictly increasing in $(0,\sqrt{2}a)$ and strictly decreasing in
$(\sqrt{2}a,2a)$. In either of these subintervals, the change of variable is defined by
$$r(s)=\int\frac{A'(s)}{P(s)}\, \dd s = \frac{A(s)}{4a}+C,$$ for a constant $C\in\R$. Fixing $C=0$, we merely have
$A[s(r)]=4ar$ and $P[s(r)]=4a$. Moreover, we can easily see that $r$ represents half of the radius of the inscribed circle (see final remark in
Section~\ref{sec:givr}).
\end{example}

\begin{remark}
The Minkowski's concept of surface area (see e.g.\ Bonnesen and Fenchel \cite[\S31]{BoFe87}), which is based on the derivative relationship
(\ref{eq:derrel}), is worth particular mention here. Let $R\in\R^d$ be a convex body of volume $V$ and surface area $A$. For any $s>0$, the
Minkowski sum
$$R(s) := R+s B^d = \{x\in\R^d\mid {\rm dist}(x,R)\leqslant s\},$$ where $B^d$ is the $d$-dimensional unit ball,
is called the {\it outer parallel body} of $R$ at distance $s$ or, equivalently, the {\it $s$-neighborhood} of $R$. According to the {\it
Steiner formula} (see e.g.\ Leichtwei\ss\ \cite[p.~30]{Le98} and Schneider \cite[Chapter~4]{Sc93}), its volume can be expressed as a polynomial
of degree $d$ in $s$, namely
$$V(s)=\sum_{i=0}^d s^i \kappa_i\, V_{d-i}\, ,$$ where
$\kappa_i$ is the volume of the $i$-dimensional unit ball, with $\kappa_0=1$, and $V_{d-i}$ is the {\it intrinsic $(d-i)$-volume of $R$}, with
special cases $V_{d}=V$ (volume of $R$) and $\kappa_1V_{d-1}=A$ (area of $R$). It is then clear that
$$\lim_{s\to +0}\frac{V(s)-V}s=\lim_{s\to +0} \frac{\dd V(s)}{\dd s} = A$$ and hence (see also Guggenheimer \cite[Chapter 4]{Gu77})
$$\frac{\dd V(s)}{\dd s}=\lim_{t\to +0}\frac{\dd V(s+t)}{\dd t}=
A(s)$$ since $V(s+t)$ is the volume of the sum $R(s)+t B^d$. We therefore retrieve Eq.~(\ref{eq:derrel}) with $r=s$.
\end{remark}

\section{Homogeneous families}

Let ${\cal F}$ be a one-parameter smooth family of compact regions in $\R^d$. Assume that $E=\R_+$ and that the parameter $s$ represents a
linear dimension of region $R(s)$, e.g., a diameter or an edge length. Then, under a dilation $s\mapsto ts$, the volume and area of that region
are clearly magnified by the factors $t^d$ and $t^{d-1}$, respectively. This means that the functions $V(s)$ and $A(s)$ fulfill the functional
equations
$$V(ts)=t^dV(s)\quad\mbox{and}\quad A(ts)=t^{d-1}A(s) \qquad (s,t\in\R_+),$$ and hence are homogeneous functions
of degrees $d$ and $d-1$, respectively, i.e., of the form
$$
V(s) = k_1 s^d \quad\mbox{and}\quad A(s) = k_2 s^{d-1} \qquad (s\in \R_+),
$$
where $k_1$ and $k_2$ are positive constants.

Starting from this observation, Tong \cite{To97} noted that, for such homogeneous functions, the derivative relationship (\ref{eq:derrel}) holds
for the change of variable
\begin{equation}\label{eq:Tongc}
r(s)=d\, \frac{V(s)}{A(s)}
\end{equation}
and the new variable $r$ also represents a linear dimension.

Note, however, that formula (\ref{eq:Tongc}) can also be valid for families of non-similar regions (see Example \ref{ex:hex}).

\begin{definition}
A smooth family ${\cal F}$ is said to be {\it homogeneous} if the change of variable in (\ref{eq:Tongc}) ensures relation (\ref{eq:derrel}).
This change of variable is then called the {\em inradius function} of ${\cal F}$.
\end{definition}


The following proposition yields equivalent conditions for a smooth family to be homogeneous.

\begin{proposition}\label{prop:tong}
Let ${\cal F}$ be a smooth family of compact regions in $\R^d$ and let $r(s)$ be given by Eq.~(\ref{eq:r}). Then the following assertions are
equivalent:
\begin{description}
\item{$i)$} There exists a constant $C\in\R$ such that
$$
r(s)=d\, \frac{V(s)}{A(s)}+C \qquad (s\in E).
$$
\item{$ii)$} There exists a constant $k>0$ such that
$$
A(s)^d=k V(s)^{d-1} \qquad (s\in E).
$$
\item{$iii)$} There exists a differentiable change of variable $\phi:E\to\phi(E)$ and constants $k_1,k_2>0$ such that
$$
V(s)=k_1\phi(s)^d \quad \mbox{and} \quad A(s)=k_2\phi(s)^{d-1}\qquad (s\in E).
$$
\end{description}
\end{proposition}

\begin{proof}
$i)\Leftrightarrow ii)$ Since $V(s)$ is differentiable, so is $A(s)$. Then, from Eq.\ (\ref{eq:rp}), we have
\begin{eqnarray*}
\exists\, C\in\R\, :\, r(s)=d\, \frac{V(s)}{A(s)}+C
& \Leftrightarrow & d\,\frac{\dd}{\dd s}\,\frac{V(s)}{A(s)}=\frac{V'(s)}{A(s)}\\
& \Leftrightarrow & d\,\frac{\dd}{\dd s}\ln A(s)=(d-1)\,\frac{\dd}{\dd s}\ln V(s)\\
& \Leftrightarrow & \exists\, k>0\, :\, A(s)^d=k V(s)^{d-1}.
\end{eqnarray*}

$ii)\Rightarrow iii)$ For any $s\in E$, define $\phi(s)=V(s)^{1/d}$. Then $V(s)=\phi(s)^d$ and $$A(s)= k^{1/d}\, V(s)^{\frac{d-1}d} =
k^{1/d}\,\phi(s)^{d-1}.$$

$iii)\Rightarrow ii)$ Clear.
\end{proof}

According to assertion $(ii)$, Eq.~(\ref{eq:Tongc}) forces the functions $A(s)^d$ and $V(s)^{d-1}$ to be linearly dependent in $E$. Thus, it
turns out that a family is homogeneous if and only if the {\it isoperimetric ratio}
\begin{equation}\label{eq:ir}
Q(s)=\frac{A(s)^d}{V(s)^{d-1}}
\end{equation}
(introduced in P\'olya \cite{Po54}) is a constant function on $E$.

On the other hand, assertion $(iii)$ clearly means that $V(s)$ and $A(s)$ are homogeneous functions of degrees $d$ and $d-1$, respectively, up
to the same change of variable $\phi(s)$. This justifies the terminology ``homogeneous family''. Clearly, this function $\phi(s)$ represents a
linear dimension and identifies with $V(s)^{1/d}$ up to a positive multiplicative constant.

We have seen in the beginning of this section that any smooth family of similar regions is homogeneous whenever the parameter $s$ represents a
linear dimension. The following corollary shows that this property holds even if $s$ does not represent a linear dimension.

\begin{corollary}
If the regions of a smooth family are all similar then the family is homogeneous.
\end{corollary}

\begin{proof} Since the regions are all similar, the isoperimetric ratio (\ref{eq:ir}), which does not
depend on the size (e.g., length of diameter) of $R(s)$, is a constant function on $E$.

\noindent {\it Alternative proof}. For any $s\in E$, let $\phi(s)$ be the diameter of region $R(s)$. Since the regions are all similar, the
functions $V(s)$ and $A(s)$ are constant multiples of $\phi(s)^d$ and $\phi(s)^{d-1}$,
 respectively.
\end{proof}

The following example shows that a homogeneous family need not be constructed from similar regions, even if the transformation carrying any
region into any other one is angle-preserving.

\begin{example}\label{ex:hex}
Consider a smooth family of hexagons whose inner angles all have a fixed amplitude $2\pi/3$ and the consecutive sides have lengths $a(s)$,
$b(s)$, $c(s)$, $a(s)$, $b(s)$, and $c(s)$, respectively. Then it can be easily shown that
\begin{eqnarray*}
A(s) &=& \frac{\sqrt{3}}2[a(s)b(s)+b(s)c(s)+c(s)a(s)], \\
P(s) &=& 2[a(s)+b(s)+c(s)].
\end{eqnarray*}
By choosing $a(s)=1$, $b(s)=s^2$, and $c(s)=(s+1)^2$, where $s\in\R_+$, we  see that this particular family of hexagons is homogeneous.
Moreover, even though the interior angles are fixed, the hexagons are not similar since the functions $a(s)$, $b(s)$, and $c(s)$ are not
linearly dependent.
\end{example}

Before closing this section, let us present an alternative interpretation of homogeneous families.


Introduced in economics, the concept of {\it elasticity} is defined as the proportional (or percent) change in one variable relative to the
proportional change in another variable. For example, the {\it price elasticity of demand} measures the change in quantity demanded with respect
to the change in price.

Applying this concept to the volume function $V(s)$ and the inradius function $r(s)$, defined in Eq.~(\ref{eq:Tongc}),
%
%
%
we can define the $r$-elasticity of volume as the proportional change in volume relative to the proportional change in the linear dimension $r$,
that is, in view of Eq.~(\ref{eq:rp}),
$$e_{V,r}(s):= \frac{~\frac{\dd V(s)}{V(s)}~}{\frac{\dd r(s)}{r(s)}} = \frac{V'(s)}{r'(s)}\, \frac {r(s)}{V(s)}= \frac {r(s)\, A(s)}{V(s)}$$
and we observe immediately that a smooth family is homogeneous if and only if
$$e_{V,r}(s)=d.$$

Considering the family of rhombi in Example~\ref{ex:vstrmon}, we simply have $e_{A,r}(s)=1$, which shows that the elasticity may be constant
while being different from $d$. Notice also that a unit elasticity means that if $r$ increases by $x$ percent then so does the area.

\section{Finding homogeneous families}

Consider an $n$-parameter family of compact regions in $\R^d$ with boundaries of finite measures,
$$
{\cal C}:=\{R(x) \subset \R^d\mid x\in F\},
$$
where $F:=F_1\times\cdots\times F_n$ is the product of $n$ open intervals of the real line. We assume that with this family is associated a
differentiable function $V:F\to\R_+$ and a continuous function $A:F\to\R_+$ such that, for any $x\in F$, the values $V(x)$ and $A(x)$ represent
respectively the volume and the surface area of region $R(x)$.

We will call such a family an $n$-parameter {\em smooth}\/ family.

\begin{example}\label{ex:par}
The class of all parallelograms in $\R^2$ can be regarded as a 3-parameter smooth family of compact figures, which can be parameterized by side
lengths $x_1>0$ and $x_2>0$, and an angle $x_3\in (0,\pi)$. In this case the corresponding area and perimeter functions are respectively given
by
$$
A(x)=x_1x_2\sin x_3 \quad \mbox{and} \quad P(x)=2x_1+2x_2.
$$
\end{example}

In this section we investigate the following problem. Given an $n$-parameter family ${\cal C}$ as defined above, find homogeneous subfamilies,
if any.

More formally, we are searching for differentiable curves
\begin{equation}\label{eq:curve}
x:E\to F,
\end{equation}
with an appropriate open real interval $E$, such that the one-parameter family
$$
\{R[x(s)] \mid s\in E\}
$$
is smooth and homogeneous.

Clearly, smoothness is ensured as soon as the function $V[x(s)]$ is strictly monotone in $E$. According to Proposition~\ref{prop:tong},
homogeneity is ensured if (see Eq.~(\ref{eq:ir}))
$$
Q[x(s)]=k \qquad (s\in E)
$$
for some $k>0$, where $Q(x):=A(x)^d/V(x)^{d-1}$. This means that the equation
$$
Q(x)=k \qquad (x\in F)
$$
represents a level hypersurface in $F$ and each differentiable curve (\ref{eq:curve}) along that hypersurface represents a homogeneous family
associated with the constant $k$.

The admissible values of $k$ are given by the well-known $d$-dimensional {\it isoperimetric inequality} (see e.g.~\cite{Ba80}, \cite{BuZa88},
\cite{Ch01}), which states that if $R$ is a compact domain in $\R^d$ with piecewise smooth boundary then
\begin{equation}\label{eq:isoineq}
\frac{A^d}{V^{d-1}} \geqslant d^d \kappa_d
\end{equation}
where $V$ and $A$ are respectively the volume and the area of $R$ and
$$\kappa_d:=\frac{\pi^{d/2}}{\Gamma(d/2+1)}$$ is the volume of the $d$-dimensional
unit ball. Here the equality sign in (\ref{eq:isoineq}) holds if and only if $R$ is the $d$-dimensional unit ball.

Thus, the constant $k$ is bounded below by $d^d \kappa_d$. For example, for $d=2$ and $d=3$, this lower bound is given by $4\pi$ and $36\pi$,
respectively.


For a particular $n$-parameter smooth family ${\cal C}$ of regions in $\R^d$, we have to refine the lower bound of constant $k$ by calculating
$$k_{\min}({\cal C}) = \inf_{x\in F} Q(x)$$
which, of course, does not depend on the parameterization of family ${\cal C}$.

For example, if ${\cal C}$ is the class of all $n$-gons in $\R^2$, we have $$k_{\min}({\cal C}) = 4n\tan(\pi/n)
$$ and this bound is achieved for the regular $n$-gons. In other words, the isoperimetric inequality for
$n$-gons $R$ in $\R^2$ is $$\frac{P^2}{A} \geqslant 4n\tan(\pi/n)$$ with equality if and only if $R$ is regular. In Table~\ref{tab:t1} we list
results for some other examples. See also Florian \cite{Fl93} and Mitrinovi\'c {\it et al.}\ \cite[Chapter~20]{MiPeVo89} for recent surveys on
isoperimetric inequalities for polytopes.

\begin{table}[tbp]\centering
\begin{tabular}{|c|l|l|c|}
\hline
$d$ & Class ${\cal C}$ & Optimal regions & $k_{\min}({\cal C})$ \\
\hline
2 & triangles & equilateral triangles & $12\sqrt{3}^{\mathstrut}$ \\
2 & right triangles & isosceles triangles & $2(2+\sqrt{2})^2$ \\
2 & $n$-gons & regular $n$-gons & $4n\tan(\pi/n)$ \\
3 & rectangular parallelepipeds & cubes & $216$ \\
3 & right circular cylinders & height = diameter & $54\pi$ \\
3 & right circular cones & height = $\sqrt{2}$ $\times$ diameter & $72\pi$ \\
3 & right square pyramids & height = $\sqrt{2}$ $\times$ side & 288 \\
3 & regular tori & apples with $r_1=r_2$ & $16\pi^2$ \\
\hline
\end{tabular}
\caption{Isoperimetric ratios for various regions} \label{tab:t1}
\end{table}

\begin{example}
Coming back to Example~\ref{ex:par} with the class ${\cal C}$ of all parallelogram in $\R^2$, we have $k_{\min}({\cal C}) = 16$. For example, if
we fix $k=32$, $x_1(s)=\sqrt{s}$, and $x_2(s)=s-\sqrt{s}$, then the third function $x_3(s)$ must be given by
$$x_3(s)=\arcsin\frac{s/8}{\sqrt{s}-1}$$
throughout the open interval $E=(24-16\sqrt{2},24+16\sqrt{2})$. Interestingly, we observe that in this case
$$A[x(s)]=s^2/8 \quad\mbox{and}\quad P[x(s)]=2s$$ are homogeneous functions.
\end{example}

Let us now investigate an interesting case. Let $m$ be an integer such that $1\leqslant m\leqslant n$ and suppose that the volume and area
functions associated with the family ${\cal C}$ are homogeneous of degrees $d$ and $d-1$ in the first $m$ variables, i.e., they fulfill the
functional equations
\begin{eqnarray}
V(tx_1,\ldots,tx_m,x_{m+1},\ldots,x_n) &=& t^d V(x_1,\ldots,x_m,x_{m+1},\ldots,x_n)
\label{eq:vh}\\
A(tx_1,\ldots,tx_m,x_{m+1},\ldots,x_n) &=& t^{d-1} A(x_1,\ldots,x_m,x_{m+1},\ldots,x_n)\label{eq:ah}
\end{eqnarray}
for all $t\in \R_+$ and all $x\in F$, where $F_1=\cdots =F_m=\R_+$. For example, the first $m$ variables $x_1,\ldots,x_m$ might represent linear
dimensions and the remaining variables $x_{m+1},\ldots,x_n$ might represent angles.

Note that these functions are necessarily of the form (see e.g.\ Acz\'el and Dhombres \cite[Chapter~20]{AcDh89})
\begin{eqnarray*}
V(x_1,\ldots,x_m,x_{m+1},\ldots,x_n) &=& x_1^d\,
f\Big(\frac{x_2}{x_1},\ldots,\frac{x_m}{x_1},x_{m+1},\ldots,x_n\Big)\\
A(x_1,\ldots,x_m,x_{m+1},\ldots,x_n) &=& x_1^{d-1}\, g\Big(\frac{x_2}{x_1},\ldots,\frac{x_m}{x_1},x_{m+1},\ldots,x_n\Big)
\end{eqnarray*}
where $f:F_2\times\cdots\times F_n\to \R_+$ and $g:F_2\times\cdots\times F_n\to \R_+$ are arbitrary continuous functions (constants if $n=1$).

Now, by using Eqs.~(\ref{eq:vh}) and (\ref{eq:ah}) with $t=1/x_1$, we immediately see that the homogeneity condition
$$
Q[x(s)]=k \qquad (s\in E),
$$
which must hold for some $k\geqslant k_{\min}({\cal C})$, is equivalent to the condition
\begin{equation}\label{eq:qzs}
Q[z(s)]=k \qquad (s\in E),
\end{equation}
where
$$
z_i(s):=\cases{\displaystyle{\frac{x_i(s)}{x_1(s)}}\, , & if $i\leqslant m$,\cr \cr x_i(s), & else.}
$$

%
%


Thus the homogeneity condition is ensured whenever we can find a differentiable curve $z:E\to F$, with $z_1(s)=1$, fulfilling (\ref{eq:qzs}).
Note that when $n=2$, (\ref{eq:qzs}) becomes $$Q[1,z_2(s)]=k \qquad (s\in E)$$ and, if the left hand side is constant in no open subinterval of
$E$ then $z_2(s)$ generally takes on a finite number of possible values.

\begin{example}
Consider the 2-parameter class of rectangles in $\R^2$. They can be parameterized, e.g., either by the length $x_1\in\R_+$ and the width
$x_2\in\R_+$, or by the half-diagonal $x_1\in\R_+$ and the angle between the diagonals $x_2\in (0,\pi/2)$. In either case, the homogeneity
condition leads to considering only similar rectangles.
\end{example}

\begin{remark}
As the searching of homogeneous subfamilies is based only on the volume and area functions, they do not depend on the parameterization used to
describe the $n$-parameter family.
\end{remark}

\section{Geometric interpretations of the inradius function}
\label{sec:givr}

As we already observed in Example~\ref{ex:rectfixw}, a geometric meaning of the variable of differentiation $r$ is not always apparent. However,
for some homogeneous families, where $r(s)$ is the inradius function, given by Eq.~(\ref{eq:Tongc}), interpretations can be found.

For example, Emert and Nelson \cite{EmNe97} proved that, for any family of similar circumscribing polytopes, the variable $r$ represents the
radius of the inscribed sphere, that is, the inradius. For an earlier work on regular polytopes, see Miller \cite{Mi78}.

Interestingly, Cohen \cite{Co65} showed that, for any $d$-dimensional circumscribing polytope of inradius $r$ and area $A$, its volume is given
by
$$
V=\frac rd\, A
$$
which corresponds to Tong formula for similar circumscribing polytopes.

Other examples have been discussed recently by Dorff and Hall \cite{DoHa}. Among these, we have the following remarkable result, that was shown
for families of similar regions in $\R^2$ and $\R^3$. We state this result in $\R^d$ and for homogeneous families. Also, Eq.~(\ref{eq:star}) was
previously unknown.

\begin{proposition}\label{prop:star}
Let ${\cal F}$ be a homogeneous family of $n$-faced polyhedra $R(s)$ that are star-like with respect to a point $T(s)$ in the interior of
$R(s)$. Let $P_i(s)$ be the pyramid whose base is the $i$th facet of $R(s)$ and whose vertex is $T(s)$. Then
\begin{equation}\label{eq:star}
r(s) = \sum_{i=1}^n \frac{A_i(s)}{A(s)}\, r_i(s)
\end{equation}
and
\begin{equation}\label{eq:star2}
\frac 1{r(s)} = \sum_{i=1}^n \frac{V_i(s)}{V(s)}\, \frac 1{r_i(s)}
\end{equation}
where $V_i(s)$, $A_i(s)$, and $r_i(s)$ are respectively the volume of $P_i(s)$, the surface area of the base of $P_i(s)$, and the altitude from
$T(s)$ of $P_i(s)$.
\end{proposition}

\begin{proof}
Since $P_i(s)$ is a $d$-dimensional pyramid, we have
$$V_i(s)=\frac 1d\, A_i(s)r_i(s).$$ Then, as the family is homogeneous, we have, by
(\ref{eq:Tongc}), $$r(s)=d\,\frac{V(s)}{A(s)}=d\,\sum_{i=1}^n \frac{V_i(s)}{A(s)} = \sum_{i=1}^n \frac{A_i(s)}{A(s)}\, r_i(s),$$ which proves
(\ref{eq:star}) and $$\frac 1{r(s)}=\frac 1d\,\frac{A(s)}{V(s)} = \frac 1d\,\sum_{i=1}^n \frac{A_i(s)}{V(s)} = \sum_{i=1}^n
\frac{V_i(s)}{V(s)}\, \frac 1{r_i(s)},$$ which proves (\ref{eq:star2}).
\end{proof}

Eq.~(\ref{eq:star}) simply means that the variable of differentiation $r(s)$ is the arithmetic mean of the altitudes from $T(s)$ of the pyramids
$P_i(s)$, weighted by the relative areas of the corresponding facets. Similarly, Eq.~(\ref{eq:star2}) means that $r(s)$ is the harmonic mean of
the altitudes from $T(s)$ of the pyramids $P_i(s)$, weighted by the relative volumes of these pyramids. Particularly, these both means do not
depend upon the choice of $T(s)$.

Clearly, Proposition~\ref{prop:star} generalizes Emert and Nelson's result and gives an interpretation of the Tong inradius (\ref{eq:Tongc}) as
an average inradius for non-circumscribing star-like regions. In some sense this justifies the terminology ``inradius function''.


For convex polytopes, Eq.~(\ref{eq:star}) can be generalized as follows. Let $R\subseteq\R^d$ be an $n$-faced convex polytope and let
$h_R:\R^d\to \R$ be its {\it support function}: $$h_R(u)=\max \{x\cdot u\mid x\in R\},$$ where $\cdot$ denotes the standard inner product on
$\R^d$. Then, assuming that $R$ has facet unit normals $u_1,\ldots,u_n$ and corresponding facet areas $A_1,\ldots,A_n$, we have (see e.g.\
Leichtwei\ss\ \cite[p.~22]{Le98})
$$V_R=\frac 1d\sum_{i=1}^n A_i \, h_R(u_i).
$$
Considering a homogeneous family of polytopes, we have immediately $$r(s)=\sum_{i=1}^n \frac{A_i(s)}{A(s)}\, h_{R(s)}[u_i(s)],$$ showing that
$r(s)$ is the weighted arithmetic mean of the functions $h_{R(s)}[u_i(s)]$.

Notice that any compact star-like set $K$ can be approximated arbitrarily closely by star-like polyhedra $\{P_i\}$ so that the volume and
surface of the $P_i$ tend in the limit to the volume and area of $K$. So the above results of this section can be easily extended to compact
star-like sets.

The harmonic mean is also encountered when considering a right cylinder in $\R^d$ obtained by appropriately lifting a region embedded in
$\R^{d-1}$.

\begin{proposition}
Let ${\cal F}_{d-1}$ be a homogeneous family of compact regions in $\R^{d-1}$ with inradius function $r_{d-1}(s)$. Consider the homogeneous
family ${\cal F}_d$ of right cylinders in $\R^d$ obtained by orthogonally lifting each region of ${\cal F}_{d-1}$ to a height of $2r(s)$. Then
the inradius function of ${\cal F}_d$ is given by
$$
r_d(s)=H_d[r_{d-1}(s),\ldots,r_{d-1}(s),r(s)],
$$
where $H_d$ is the $d$-variable symmetric harmonic mean.
\end{proposition}

\begin{proof}
The result immediately follows from the equalities
\begin{eqnarray*}
V_{d}(s) &=& 2V_{d-1}(s)r(s) \\
A_{d}(s) &=& 2V_{d-1}(s)+2A_{d-1}(s)r(s) = 2V_{d-1}(s)+2(d-1)\frac{V_{d-1}(s)}{r_{d-1}(s)}\, r(s),
\end{eqnarray*}
where $V_{d-1}(s)$, $A_{d-1}(s)$, $V_{d}(s)$, and $A_{d}(s)$ denote the volume and area functions of ${\cal F}_{d-1}$ and ${\cal F}_{d}$,
respectively.
\end{proof}

\begin{remark}
It can be easily proved that the function $A(s)$ of a one-parameter smooth family ${\cal F}$ is constant if and only if
$$r(s)=\frac{V(s)}{A(s)}+C.$$ When $C=0$, the function $r(s)$ identifies with one
$d$th of the inradius function and the elasticity $e_{V,r}(s)$ is one. If, moreover, the regions are star-like $n$-faced polyhedra as in
Proposition~\ref{prop:star}, then Eqs.~(\ref{eq:star}) and (\ref{eq:star2}) become respectively
$$r(s) = \frac 1d\, \sum_{i=1}^n \frac{A_i(s)}{A(s)}\, r_i(s)$$ and $$\frac 1{r(s)} = d\,\sum_{i=1}^n
\frac{V_i(s)}{V(s)}\, \frac 1{r_i(s)}$$ Example~\ref{ex:vstrmon} illustrates these latter two formulas.
\end{remark}

\section{Bonnesen-style inequalities with Tong inradius}

Let $R$ be any compact plane figure in $\R^2$ with piecewise smooth boundary. Denote by $P$ and $A$ its perimeter and area, respectively. The
isoperimetric inequality (\ref{eq:isoineq}) ensures that the quantity
$$P^2-4\pi A,$$ known as the {\em isoperimetric deficit} of $R$, is nonnegative.
Bonnesen (\cite{Bo21},\cite{Bo26},\cite{Bo29},\cite{Fu91}) found lower bounds for the isoperimetric deficit by establishing the following
inequalities:
\begin{eqnarray*}
P^2-4\pi A & \geqslant & (P-2\pi r)^2, \\
P^2-4\pi A & \geqslant & (\frac{A}r-\pi r)^2,
\end{eqnarray*}
where $r$ is the radius of any circle inscribed in $R$. Interestingly, Osserman \cite{Os79} proved that these inequalities are each
algebraically equivalent to
$$rP \geqslant A+\pi r^2.$$

It is simple routine that all these inequalities still hold when $r$ is the Tong inradius $r=2\frac{A}{P}$.

For a general dimension $d$, we have the following result.

\begin{proposition}
Let $R$ be a compact domain in $\R^p$ with piecewise smooth boundary. Denote by $A$ and $V$ its area and volume, respectively, and let
$r=d\,\frac{V}{A}$ be its Tong inradius. Then we have
\begin{eqnarray}
A^d-d^d\kappa_d V^{d-1} & \geqslant & (A-d\,\kappa_d\, r^{d-1})^d, \label{eq:bon1}\\
A^d-d^d\kappa_d V^{d-1} & \geqslant & \Big(\frac{V}{r}-\kappa_d\, r^{d-1}\Big)^d, \label{eq:bon2}
\end{eqnarray}
and
\begin{equation}\label{eq:oss}
r A \geqslant V+(d-1)\, \kappa_d\, r^d.
\end{equation}
\end{proposition}

\begin{proof}
Inequality (\ref{eq:bon1}) is immediate if we observe that the right hand side writes
$$\Big(A-d^d\,\kappa_d\, \frac{V^{d-1}}{A^{d-1}}\Big)^d = \frac 1{A^{d(d-1)}}\Big(A^d-d^d\,\kappa_d\, V^{d-1}\Big)^d.$$
The other inequalities are routine.
\end{proof}

\section{Conclusion}
We have explored the idea of the derivative of the volume of a region in $\R^d$ with respect to some variable $r$ equaling its surface area for
homogeneous families. This area of investigation is intriguing and appears not to have been previously studied. We have just skimmed the
surface, and there are a lot of questions to be answered. For example, what other geometric interpretations are there for the inradius function?

\section*{Acknowledgements}

The authors are grateful to Jean-Paul Doignon for calling their attention to the Minkowski's concept of surface area.

\noindent
Jean-Luc Marichal\\
Applied Mathematics Unit, University of Luxembourg\\
162A, avenue de la Fa\"{\i}encerie, L-1511 Luxembourg, Luxembourg\\
Email: jean-luc.marichal[at]uni.lu\\

\noindent
Michael Dorff\\
Department of Mathematics, 281 TMCB, Brigham Young University\\
Provo, Utah 84602, U.S.A.\\
Email: mdorff[at]math.byu.edu

\begin{thebibliography}{99}

\bibitem{AcDh89} J.\ Acz\'el and J.\ Dhombres, {\em Functional equations in several variables with applications to mathematics, information theory and to the natural and social sciences} (Cambridge Univ.\ Press, Cambridge, 1989).


\bibitem{Ba80} C.\ Bandle, {\it Isoperimetric inequalities and applications} (Pitman, Boston, 1980).

\bibitem{Bo21} T.\ Bonnesen, Uber eine versch\"arfung der isoperimetrischen
ungleichheit des kreises in der ebene und auf der kugeloberfl\"ache nebst einer anwendung auf eine Minkowskische ungleichheit f\"ur konvexe
k\"orper, {\em Math. Annalen} {\bf 84} (1921) 216--227.

\bibitem{Bo26} T.\ Bonnesen, Quelques probl\`emes isop\'erim\'etriques, {\em Acta
Mathematica} {\bf 48} (1926) 123--178.

\bibitem{Bo29} T.\ Bonnesen, {\em Les probl\`emes des isop\'erim\`etres et des
is\'epiphanes} (Gauthier-Villars, Paris, 1929).

\bibitem{BoFe87} T.\ Bonnesen and W.\ Fenchel, {\it Theory of convex bodies} (BCS Associates, Moscow, Idaho, 1987).

\bibitem{BuZa88} Y.D.\ Burago and V.A.\ Zalgaller, {\it Geometric inequalities} (Springer-Verlag, Berlin, 1988).

\bibitem{Ch01} I.\ Chavel, {\it Isoperimetric inequalities} (Cambridge University Press, Cambridge, 2001).

\bibitem{Co65} M.J.\ Cohen, Solution of elementary problems: E1671 -- Ratio of volume of inscribed sphere to polyhedron, {\it Amer.\ Math.\ Monthly} {\bf 72} (2) (1965) 183--184.

\bibitem{DoHa} M.\ Dorff and L.\ Hall, Solids in $\mathbf{R}^n$ whose area is the derivative of the volume, {\it College Math. J.} {\bf 34} (5) (2003) 350--358.

\bibitem{EmNe97} J.\ Emert and R.\ Nelson, Volume and surface area for polyhedra and polytopes, {\it Math.\ Mag.} {\bf 70} (5) (1997) 365--371.

\bibitem{Fl93} A.\ Florian, Extremum problems for convex discs and polyhedra, {\it Handbook of convex geometry}, Vol.\ A, B (North-Holland, Amsterdam, 1993) 177--221.

\bibitem{Fu91} B.\ Fuglede, Bonnesen's inequality for the isoperimetric deficiency of closed curves in the plane, {\it Geometriae Dedicata} {\bf 38} (3) (1991) 283--300.

\bibitem{Gu77} H.\ Guggenheimer, {\it Applicable geometry} (Robert E.\ Krieger Publishing Co., Inc., Huntington, N.~Y., 1977).

\bibitem{Le98} K.\ Leichtwei\ss, {\it Affine geometry of convex bodies} (Johann Ambrosius Barth Verlag, Heidelberg, 1998).

\bibitem{Mi78} J.I.\ Miller, Differentiating area and volume, {\it Two-Year College Math.\ J.} {\bf 9} (1) (1978) 47--48.

\bibitem{MiPeVo89} D.S.\ Mitrinovi\'c, J.E.\ Pe\v{c}ari\'c, and V.\ Volenec, {\it Recent Advances in Geometric Inequalities} (Kluwer, Dordrecht, 1989).

\bibitem{Os79} R. Osserman, Bonnesen-style isopermetric inequalities, {\em Amer.
Math. Monthly} {\bf 86} (1979) 1--29.

\bibitem{Po54} G.\ P\'olya, {\it Induction and analogy in mathematics} (Princeton University Press, Princeton, 1954).

\bibitem{Sc93} R.\ Schneider, {\it Convex bodies: the Brunn-Minkowski theory} (Cambridge University Press, Cambridge, 1993).

\bibitem{St90} K.A. Struss, Exploring the volume-surface area relationship, {\it College Math.\ J.} {\bf 21} (1) (1990) 40--43.

\bibitem{To97} J.\ Tong, Area and perimeter, volume and surface area, {\it College Math.\ J.} {\bf 28} (1) (1997) 57.

\end{thebibliography}
\end{document}